%%%%%%%%%%%%%%%%%%%%%%%%%%%%%%%%%%%%%%%%%%%%%%%%%%%%%%%%%%%%%%%%%
%  Version 1.0 started on 10/2/01, finished 
%  Version 1.3 revised on 3/9/03 for submission to Math.Ann.
%  Version 2.0 revised on 11/12/03 according to referee's comments
%%%%%%%%%%%%%%%%%%%%%%%%%%%%%%%%%%%%%%%%%%%%%%%%%%%%%%%%%%%%%%%%%

\input amstex
\documentstyle{amsppt}

\hcorrection{19mm}

\nologo
\NoBlackBoxes

%\NoRunningHeads
%\mag=1200
%\hsize=31 pc
%\vsize=44 pc
%\baselineskip 15pt
%\hcorrection{5mm}

\topmatter 

\title The Classification of Dehn fillings on the outer torus of a 1-bridge braid
exterior which produce solid tori
\endtitle
\author  Ying-Qing Wu$^1$
\endauthor
\leftheadtext{Ying-Qing Wu}
\rightheadtext{Dehn fillings on 1-bridge braid exterior}
\address Department of Mathematics, University of Iowa, Iowa City, IA
52242
\endaddress
\email  wu\@math.uiowa.edu
\endemail

\thanks  $^1$ Partially supported by NSF grant \#DMS 0203394
\endthanks

\abstract Let $K= K(w,b,t)$ be a 1-bridge braid in a solid torus $V$,
and let $\gamma$ be a $(p,q)$ curve on the torus $T = \partial V$ of
the exterior $M_K$ of $K$.  It will be shown that Dehn filling on $T$
along $\gamma$ produces a solid torus if and only if $p$ and $q$
satisfy one of four conditions determined by the parameters $(w,b,t)$
of the knot $K$.  This solves the classification problem raised by
Menasco and Zhang for such Dehn fillings.  \endabstract

\endtopmatter

\document

\define\proof{\demo{Proof}}
\define\endproof{\qed \enddemo}
\define\a{\alpha}
\redefine\b{\beta}

\define\r{\gamma}
\redefine\e{\epsilon}
\redefine\bdd{\partial}
\define\Int{\text{\rm Int}}
\input epsf.tex

\head 1. Introduction \endhead

A knot $K$ in a 3-manifold $M$ is a 0-bridge knot if it is isotopic to
a simple closed curve on $\bdd M$, and it is a {\it 1-bridge knot\/}
if it is not 0-bridge, and is isotopic to a curve $\a \cup \b$, where
$\a\subset \bdd M$, and $\b$ is a trivial arc in $M$ in the sense that
it is rel $\bdd$ isotopic to an arc on $\bdd M$.  These knots have
been related to many examples of exceptional Dehn surgery on
hyperbolic knots, and have been studied quite extensively, see for
example [Ga1, Ga2, Be, Eu, MZ, Wu1].  In particular, it was proved
by Gabai [Ga1, Ga2] that if some surgery on a knot in a solid torus
yields a solid torus then the knot must be a 1-bridge braid.  Berge
completely classified all such surgeries in [Be].  Menasco and Zhang
[MZ] studied Dehn fillings on the outer torus of 1-bridge braids,
showing that some of those Dehn fillings produce solid tori.  The main
purpose of this paper is to solve a problem raised in their paper [MZ,
Problem 7], which asked for a complete classification of all such Dehn
fillings.

Let $B_w$ be the braid group on $w$ strands, and let $\sigma_i$ be the
standard generators of $B_w$.  See [Bi] for definitions.  A braid
$\sigma$ is represented by a set of $w$ strings in $D^2 \times I$.
Thus when gluing $D^2 \times 0$ to $D^2 \times 1$, we obtain the
closure of $\sigma$, which is a knot or link in the solid torus $V =
D^2 \times S^1$.  A {\it 1-bridge braid\/} is a knot $K = K(w,b,t)$ in
$V$ which is the closure of the braid

$$\sigma(w,b,t) = \sigma_{b} \cdots \sigma_2 \sigma_1 (\sigma_{w-1}
\cdots \sigma_2 \sigma_1)^t.$$ See Figure 1.1 for the braid
$\sigma(7,4,2)$.  Note that not all 1-bridge knots in $V$ are 1-bridge
braids.  Note also that if $b=0$ or $w-1$ then $K(w,b,t)$ is isotopic
to a curve on $T =\bdd V$, hence is a 0-bridge knot, in which case the
exterior of $K$ is a cable space.  Dehn filling on cable spaces are
well understood, see [Go].  Thus we will restrict our attention to
$K(w,b,t)$ with $1\leq b \leq w-2$.  Up to homeomorphism of $V$
obtained by twisting along meridional disks, we may also assume that
$1\leq t \leq w-1$.

\bigskip
\leavevmode

\centerline{\epsfbox{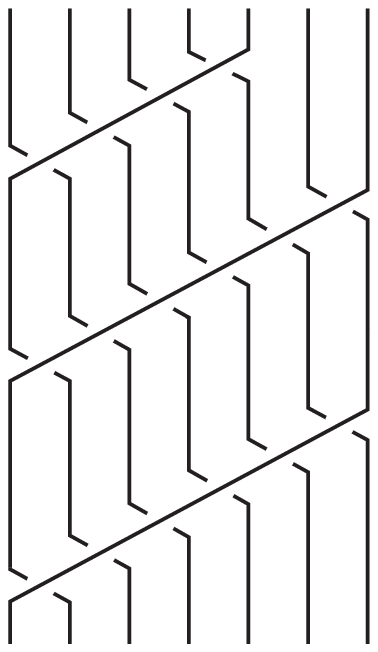}}
\bigskip
\centerline{Figure 1.1}
\bigskip

Let $\r$ be a $(p,q)$ curve on $T$ with respect to the standard
longitude-meridian pair of $V$.  Let $\a \cup \b$ be a 1-bridge
presentation of a 1-bridge knot $K$ in $V$.  Thus $\a \subset T$, and
$\b$ is a trivial arc in $V$.  Up to isotopy we may assume that $\a$
is disjoint from $\r$.  It can be shown (Lemma 2.5) that the manifold
$M_K(\r)$ obtained by performing a Dehn filling on the exterior $M_K$
of $K$ along the curve $\r$ is homeomorphic to the manifold $X[\r]$
obtained by attaching a 2-handle to the genus 2 handlebody $X= V -
\Int N(\b)$ along the curve $\r$.  Thus the problem of determining
which $M_K(\r)$ is a solid torus is the same as to determine which
$X[\r]$ is a solid torus.

A trivial arc $\b$ in a solid torus is {\it $\r$-trivial\/} if it is
isotopic rel $\bdd$ to an arc $\b'$ on $T$ intersecting $\r$ always in
the same direction, in which case the minimal intersection number
between all such $\b'$ and $\r$ is called the {\it jumping number\/}
of $\b$ with respect to $\r$.  See Definition 2.1 for more details.
The main theorem of Section 2 states that $X[\r]$ is a solid torus if
and only if $\b$ is $\r$-trivial, and the jumping number of $\b$ with
respect to $\r$ is either 1 or $\min (q, p-q)$.  See Theorem 2.2.  The
result will be used in [Wu2] to classify completely tubing
compressible tangles.

Because of the relation between $M_K(\r)$ and $X[\r]$, the above
theorem gives a necessary and sufficient condition for $M_K(\r)$ to be
a solid torus.  See Corollary 2.6.  It will also be shown (Theorem
2.8) that $M_K (\r)$ is reducible if and only if $K$ is a cable of a
knot parallel to $\r$.  We remark that in general $M(\r)$ may not be a
solid torus for any $\r$ since we are performing Dehn filling on the
outer torus.  Because of that, the powerful Reducible Surgery Theorem
of Scharlemann [Sch, Theorem 6.1] does not apply to this situation.

The results in Section 2 will be used in Section 3 to classify all
1-bridge knots in $V$ which admits a Dehn filling on the outer torus
producing a solid torus.  We will also determine all such Dehn filling
slopes for a given 1-bridge braid $K(w,b,t)$.  It will be shown that
$M_K(\r)$ is a solid torus if and only if $p$ and $q$ satisfy one of
four conditions determined by the parameters of the 1-bridge braid
$K=K(w,b,t)$.  See Theorem 3.6.  This solves the classification
problem raised by Menasco and Zhang for such Dehn fillings [MZ,
Problem 7].  Some computational results based on these theorems will
be given at the end of that section.

We remark that 1-bridge knots are not the only ones which may admit
some solid torus Dehn fillings on the outer tori of their exteriors.
Given a knot $K$ in a solid torus $V$ and a slope $\r$ on $\bdd V$,
let $K'$ be the core of the Dehn filling solid torus in $V(\r)$.  Then
$M_K(\r)$ is a solid torus if and only if the link $L = K \cup K'$ is
a generalized Brunnian link in $V(\r)$ in the sense that the
complement of each component of $L$ is a solid torus; in particular,
if $L= K \cup K'$ is a Brunnian link in $S^3$ and $V= S^3 - \Int
N(K')$ then Dehn filling along the meridian slope of $K'$ on $M_K$ is
a solid torus.  However, if $M_K$ admits two solid torus Dehn fillings
on the outer torus, then by Gabai's theorem [Ga1] the core of such a
Dehn filling solid torus is a 1-bridge braid, in which case it is easy
to show that $K$ must be a 1-bridge knot, and hence a 1-bridge braid
by Corollary 2.6.

We work in the smooth or piecewise linear category.  Denote by $N(Y)$
a closed regular neighborhood of a subset $Y$ in a 3-manifold $M$, and
by $|Y|$ the number of components of $Y$.  Given a knot $K$ in a solid
torus $V$, let $M_K = V - \Int N(K)$.  Let $\r$ be a $(p,q)$-curve on
$T$ disjoint from $\bdd \b$, i.e., $\r$ represents $p[l] + q[m]$ in
$H_1(T)$, where $l = x\times S^1$ and $m = \bdd D^2 \times y$ for some
$x\in \bdd D^2$ and $y \in S^1$.  We use both $M_K(\r)$ and $M_K(p/q)$
to denote the manifold obtained from $M_K$ by Dehn filling on $T$
along the slope $\r$.  

\head 2. A criterion \endhead

Consider a proper arc $\b$ in a solid torus $V = D^2 \times S^1$,
which is {\it trivial\/} in the sense that it is rel $\bdd \b$
isotopic to an arc $\b'$ on $T = \bdd V$.  Then $X = V -\Int N(\b)$ is
a genus 2 handlebody.  Let $\r$ be a $(p,q)$ curve on $T$.  We are
interested in the question of when the manifold $X[\r]$ obtained by
attaching a 2-handle to $X$ along $\r$ is a solid torus.  Note that
$X[\r]$ can also be obtained by removing a regular neighborhood of the
arc $\b$ from the punctured lens space $V[\r]$, i.e., $(V -\Int
N(\b))[\r] = V[\r] - \Int N(\b)$.

Since $\b$ is a trivial arc in $V$, there is a meridian disk $D$ of
$V$ containing $\b$.  Since $\r$ is a $(p,q)$-curve, there is also a
meridian disk $D'$ of $V$ such that $\bdd D'$ intersects $\r$ in
exactly $p$ points in the same direction.  However, in general one
cannot choose $D$ and $D'$ to be the same disk; in other words, it may
not be possible to find an isotopy relative to $\r$ which deforms $\b$
to an arc lying on a disk $D'$ whose boundary intersects $\r$ at $p$
points.  For example, let $K$ be an arbitrary 1-bridge braid with
1-bridge presentation $\a \cup \b$, which is not a $0$-bridge braid.
Let $\r$ be a longitude on $T$, and let $D$ be a meridian of $V$
intersecting $\r$ at a single point.  Up to isotopy we may assume $\a
\cap \r = \emptyset$.  Then $\b$ cannot be rel $\r$ isotopic to an arc
on $D$, as otherwise one can show that $\b$ would be rel $\bdd$
isotopic to an arc $\b'$ on $T$ which intersects $\r$ at a single
point, so $\a \cup \b'$ would be a simple closed curve on $T$,
contradicting the fact that $K$ is not a 0-bridge knot.

\definition{Definition 2.1} Let $\r$ be a $(p,q)$-curve on $T = \bdd
V$.  A properly embedded arc $\b$ in $V$ is {\it $\r$-trivial\/} if it
lies on a meridian disk $D$ of $V$ such that $\bdd D$ intersects $\r$
at $p$ points.  In this case $\bdd \b$ divides $\bdd D$ into two arcs
$\b'$ and $\b''$.  The smaller of the intersection numbers $|\b' \cap
\r|$ and $|\b'' \cap \r|$, is called the {\it jumping number\/} of
$\b$ with respect to $\r$, denoted by $u = u(\b, \r)$.  \enddefinition

The following theorem characterizes trivial arcs $\b$ in $V$ such that
$V[\r] - \Int N(\b)$ is a solid torus.  This should be compared with
[MZ, Proposition 3], where it does not seem to have been realized that
a trivial arc in $V$ may not be $\r$-trivial.  Because of this, the
proof to [MZ, Corollary 4] is not complete.  The following theorem
will fill the gap.

\proclaim{Theorem 2.2} Let $\b$ be a trivial arc in a solid torus $V$,
and let $\r$ be a $(p,q)$-curve on $T = \bdd V$ disjoint from $\b$.
Let $X = V - \Int N(\b)$.  Then $X[\r]$ is a solid torus if and only
if (i) $\b$ is $\r$-trivial, and (ii) the jumping number
$u(\b, \r)$ equals $1$ or $\min(q,p-q)$.  \endproclaim

We need a result which determines when the boundary of a handlebody
with a curve removed is incompressible.  The following lemma is due to
Starr [St].  An alternative proof can be found in [Wu1, Theorem 1.2
and Corollary 1.3].

\proclaim{Lemma 2.3} Let $H$ be a handlebody, and let $\r$ be a simple
closed curve on $\bdd H$.  Then $\bdd H -\r$ is incompressible if and
only if there is a set of essential disks $D_1, ..., D_k$ in $H$
cutting $\bdd H$ into a set of twice punctured disks $P_1, ..., P_r$,
such that for each $i$, we have (1) each component of $\r \cap P_i$ is
an essential arc on $P_i$, and (2) there is at least one component of
$\r \cap P_i$ connecting each pair of boundary components of $P_i$.
\endproclaim

Recall that a manifold $M$ is $\bdd$-irreducible if $\bdd M$ is
incompressible in $M$.  The following lemma proves the necessity of
(1) in Theorem 2.2.

\proclaim{Lemma 2.4} Let $\b, V, \r, X$ be as in Theorem 2.2.  If $\b$
is not $\r$-trivial, then $X[\r]$ is irreducible and
$\bdd$-irreducible.  \endproclaim

\proof Let $D_0$ be a meridian disk of $V$ containing $\b$, and let
$D_1$ be another meridian disk of $V$ disjoint from $D_0$.  We may
assume that $\r$ has been deformed by an isotopy of $V$ {\it relative
  to $\b$\/} so that $\r$ intersects $\bdd D_i$ minimally.  Let $A$ be
the annulus obtained by cutting $T = \bdd V$ along the meridian curve
$\bdd D_1$.  Then $\r \cap A$ is a set of arcs $\r_1, ..., \r_n$.  If
all $\r_i$ are essential arcs in $A$, then by an isotopy we may assume
that $\r_i$ are straight arcs from one boundary component of $A$ to
another.  But in this case $\b$ would be rel $\bdd$ isotopic to a
straight arc $\b'$ in $A$ intersecting $\r$ always in the same
direction, which would imply that $\b$ is $\r$-trivial, contradicting
the assumption.

Therefore we may assume that some component $\r'$ of $\r \cap A$ is an
inessential arc in $A$.  Since the two boundary components of $A$
contain the same number of points of $\r$, there must be another
component $\r''$ of $\r \cap A$ with both endpoints on the boundary
component of $A$ which does not contain the endpoints of $\r'$.  Let
$\Delta'$ and $\Delta''$ be the disks on $A$ cut off by $\r'$ and
$\r''$, respectively.  By the minimality of $|\r \cap \bdd D_i|$, we
see that each of $\Delta'$ and $\Delta''$ must contain exactly one
endpoint of $\b$.  Rechoosing $D_0$ if necessary, we may assume that
the intersection of $\bdd D_0$ with each of $\Delta'$ and $\Delta''$
is a single arc.  See the first figure on Figure 2.1.

\bigskip
\leavevmode

\centerline{\epsfbox{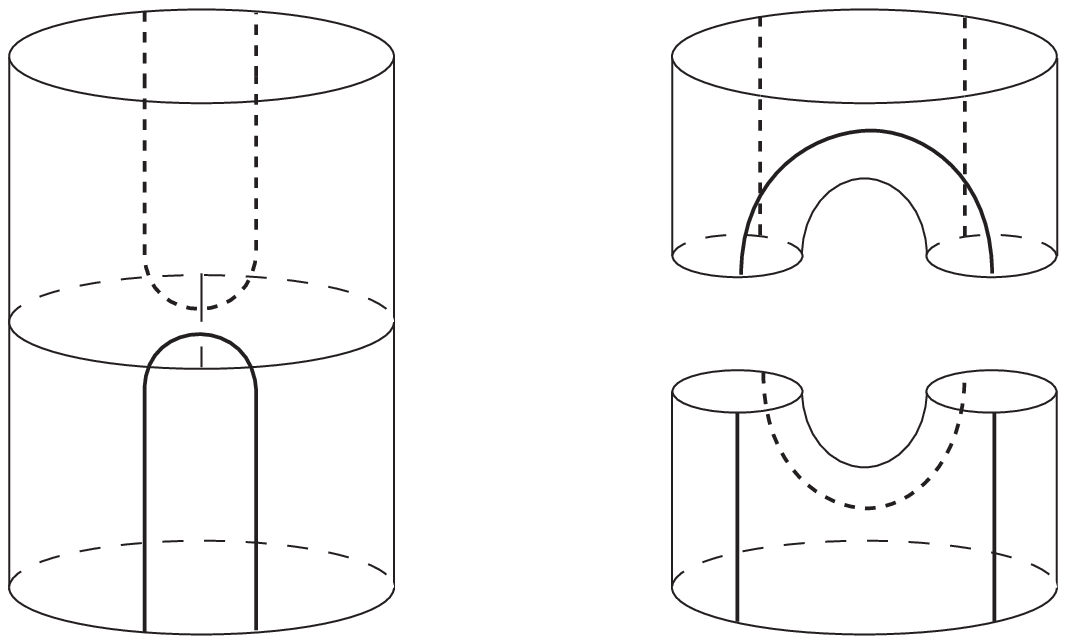}}
\bigskip
\centerline{Figure 2.1}
\bigskip

Consider the genus 2 handlebody $X = V - \Int N(\b)$.  The disk $D_0$
containing $\b$ intersects $X$ in two disks $D_2, D_3$.  Now the three
disks $D_1$, $D_2$ and $D_3$ cut the genus 2 surface $\bdd X$ into two
pairs of pants $P_1$ and $P_2$, as shown in the second figure of
Figure 2.1.  The minimality of $|\r \cap D_i|$, $i=0,1$, guarantees
that the arcs of $\r \cap P_i$ are essential.  The two arcs $\r'$ and
$\r''$ on $A$ now give rise to 6 mutually non-parallel essential
arcs, three on each of $P_1$ and $P_2$.  It follows from Lemma 2.3
that $\bdd X - \r$ is incompressible.  Since $X$ is irreducible and
$\bdd X$ is compressible, by the Handle Addition Lemma [CG], after
adding a 2-handle along $\r$, the manifold $X[\r]$ is irreducible and
$\bdd$-irreducible.  \endproof

Recall that we use $M(\r)$ to denote the Dehn filling on a torus
boundary component of $M$ along a slope $\r$.  The following lemma
relates Dehn filling to 2-handle addition.  It is true for a general
3-manifold $V$, although we will only use it for $V = D^2 \times S^1$.

\proclaim{Lemma 2.5} Let $K$ be a knot in a 3-manifold $V$ which is
isotopic to a curve $\a \cup \b$, where $\a$ is an arc on a torus
boundary component $T$ of $V$, and $\b$ is a properly embedded arc.
Let $\r$ be a simple closed curve on $T$ disjoint from $\a$.  Let $M
= V - \Int N(K)$, and $X = V - \Int N(\b)$.  Then $M(\r) = X[\r]$.
\endproclaim

\proof Consider the manifold $V(\r) = V \cup V'$, where $V'$ is the
Dehn filling solid torus.  Let $\a'$ be a proper arc in $V'$ isotopic
to $\a$ rel $\bdd \a$.  Then $K$ is isotopic to $\a' \cup \b$ in
$V(\r)$.  Let $N(D)$ be a regular neighborhood of a meridian disk $D$
in $V'$ such that $D \cap \a' = \emptyset$.  Then $V' - \Int N(D)$ is
a regular neighborhood of $\a'$ in $V(\r)$, denoted by $N(\a')$.
Hence
$$\align
M(\r) &= V(\r) - \Int N(K) \cong V(\r) - \Int(N(\b) \cup N(\a')) \\
&= (V(\r) - \Int N(\a')) - \Int N(\b) \\
&= V[\r] - \Int N(\b) = X[\r] 
\endalign 
$$
\endproof

\medskip

\demo{Proof of Theorem 2.2} Let $\a$ be an arc in $T - \r$ connecting
the endpoints of $\b$.  By Lemma 2.5 we have $X[\r] = M_K(\r) = V(\r) -
\Int N(K)$, where $K$ is a knot in $V$ isotopic to $\a\cup \b$.  Note
that $V(\r)$ is a lens space $L(p,q)$; hence by the uniqueness of
Heegaard splittings of lens spaces [BO] we see that $X[\r] = V(\r) -
\Int N(K)$ is a solid torus if and only if $K$ is a core of either $V$
or the attached solid torus $V'$ in $V(\r) = L(p,q) = V \cup V'$.

First assume that $X[\r]$ is a solid torus.  By Lemma 2.4, $\b$ is
$\r$-trivial, so it remains to show that $u = u(\b, \r) = 1$ or
$\min(q, p-q)$.  By definition $\b$ is isotopic rel $\bdd$ to an arc
$\b'$ on $T$ intersecting $\r$ transversely at $u(\b,\r)$ points in
the same direction.  Thus if we choose the core curve of the attached
solid torus $V'$ as a generator of $\pi_1 L(p,q) = \Bbb Z_p$, then the
curve $K = \a \cup \b$ represents the number $u(\b, \r)$ in $\Bbb
Z_p$.  By the above, $K$ is a core of $V'$ or $V$, which represents
$1$ or $q$ in $\pi_1 L(p,q) = \Bbb Z_p$, respectively.  Hence $u(\b,
\r) = 1$ or $\min(q,p-q)$.  This completes the proof of necessity.

Now assume that $\b$ is $\r$-trivial, and has jumping number $1$ or
$\min(q,p-q)$.  By definition there is a meridian disk $D$ of $V$
containing $\b$, with $\bdd D$ intersection $\r$ at $p$ points in the
same direction.  If $u(\b,\r) = 1$ then $\b$ is rel $\bdd$ isotopic to
an arc $\b'$ on $T$ intersecting $\r$ at a single point.  Note that in
this case $\a \cup \b'$ is a simple closed curve on $T$ and is
isotopic to a core of $V'$.  If $u(\b,\r) = \min(q,p-q)$ then since
$\r$ is a $(p,q)$ curve, there is an arc $\a'$ on $T - \r$ with $\bdd
\a' = \bdd \b$, running along the longitudinal direction of $V$ only
once in the sense that its union with an arc on $\bdd D$ is a
longitude of $V$.  Since $T-\r$ is a meridional annulus on the Dehn
filling solid torus $V'$, the two arcs $\a$ and $\a'$ are isotopic rel
$\bdd$ in $V'$; hence $K \cong \b \cup \a'$ is isotopic to a core of
$V$.  In either case $X[\r] = V(\r) - \Int N(K)$ is a solid torus.
\endproof

\proclaim{Corollary 2.6} Suppose $K$ is a 1-bridge knot in a solid
torus $V$ with 1-bridge presentation $\a \cup \b$.  Let $\r$ be a
$(p,q)$ curve on $T=\bdd V$ disjoint from $\a$.  Let $M_K = V - \Int
N(K)$.  Then $M_K(\r)$ is a solid torus if and only if (i) $\b$ is
$\r$-trivial, and (ii) the jumping number $u(\b, \r) = 1$ or $\min(q,
p-q)$.  

In particular, if $M_K(\r)$ is a solid torus then the 1-bridge knot
$K$ must be a 1-bridge braid.
\endproclaim

\proof Since $M(\r) = V[\r]$, this follows immediately from Theorem
2.2.  Note that if $K$ is isotopic to $\a \cup \b$ such that $\a$
disjoint from $\r$ and $\b$ is $\r$-trivial, then $K$ is a 1-bridge
braid.  
\endproof

\proclaim{Corollary 2.7} Let $K$ be a 1-bridge braid in a solid
torus $V$, and let $\varphi : V \to S^3$ be an embedding.  Then
$\varphi(K)$ is a nontrivial knot in $S^3$.  \endproclaim

\proof This is obvious if $\varphi(V)$ is a nontrivial torus.  So
assume $\varphi(V)$ is trivial in $S^3$, and let $\r$ be a longitude
of $V$ such that $\varphi(\r)$ bounds a meridional disk in $S^3 - \Int
\varphi(V)$.  Then $V(\r) = S^3$, so $\varphi(K)$ is trivial in $S^3$
if and only if $M_K(\r)$ is a solid torus.  By Corollary 2.6, this
implies that $K$ has a 1-bridge presentation $\a \cup \b$ such that
$\a$ is disjoint from $\r$ and $\b$ is $\r$-trivial.  Since $\r$ is a
longitude, the jumping number of $\b$ must be $0$.  It is now easy to
see that $K$ is isotopic to a simple closed curve on $T = \bdd V$, so
it is a 0-bridge knot.  Since these have been excluded from 1-bridge
braids, the result follows.  \endproof

The following theorem characterizes reducible Dehn fillings on the
outer torus of 1-bridge knots in solid tori.

\proclaim{Theorem 2.8} Let $K$ be a 1-bridge knot in a solid torus
$V$, and let $M= V - \Int N(K)$.  Let $\r$ be a $(p,q)$ curve on $T =
\bdd V$.  Then $M(\r)$ is reducible if and only if (i) $p>1$, and
(ii) $K$ is a cable of a knot $K'$ in $V$ parallel to $\r$.
\endproclaim

\proof
If $K$ is an $(r,s)$ cable of a $(p,q)$ knot in $V$, $p>1$, then $M$
is the union of a $(p,q)$-cable space $C_{p,q}$ and an $(r,s)$-cable
space $C_{r,s}$ along a boundary component.  Since $\r$ is a fiber of
the Seifert fibration of $C_{p,q}$, the Dehn filling space $C_{p,q}
(\r)$ is reducible; hence so is $M(\r) = C_{r,s} \cup C_{p,q}(\r)$.

Now assume $M(\r)$ is reducible.  Then $p>1$, for otherwise $K$ would
be a knot in $S^3$, or a knot in $S^2\times S^1$ with nontrivial
winding number, so $M(\r) = V(\r) - \Int N(K)$ would be irreducible.
This proves (i).

Now let $\a \cup \b$ be a 1-bridge presentation of $K$ such that $\a
\cap \r = \emptyset$.  By Lemma 2.5 we have $M(\r) = X[\r]$, and by
Lemma 2.4 the arc $\b$ is $\r$-trivial in $V$.  Since $M(\r) = V(\r) -
\Int N(K)$ and $V(\r)$ is a lens space, $M(\r)$ is reducible if and
only if $K$ is a knot in a ball $B$ in $V(\r)$.  Since $K$ represents
$u(\b, \r)$ in $\pi_1(V(\r))$, the jumping number $u(\b, \r) = 0$,
which implies that $K$ is isotopic to a cable of a knot $K'$ parallel
to $\r$.  
\endproof

\head 3. The classification \endhead

Let $\a \cup \b$ be a 1-bridge presentation of a 1-bridge braid $K =
K(w,b,t)$ in $V$.  Orient $K$ so that the winding number of $K$ in $V$
is $w>0$.  This induces an orientation on $\a$ and $\b$.  Cutting the
solid torus $V = D^2 \times S^1$ along a meridian disk $D$ which is
disjoint from $\b$ and intersects $\a$ minimally, the torus $T = \bdd
V$ becomes an annulus $A$.  Choose an arc $C$ on $A$ which is disjoint
from $\a$, such that $A$ cut along $C$ is a rectangle $R$ as shown in
Figure 3.1.  The intersection of $\a$ with $R$ is a set of arcs $\a_0,
..., \a_w$, ordered from left to right on $R$.  Let $\a(0)$ and
$\a(1)$ be the initial and ending points of $\a$.  We may assume that
the orientation of $\a$ points downward, and the arc $C$ has been
chosen so that the arc $\a'$ among the $\a_i$ containing $\a(0)$ is to
the left of the arc $\a''$ containing $\a(1)$.  See Figure 3.1 for an
example, where the $\a_i$ are drawn as vertical lines, and $(w,b,t) =
(7,4,2)$.

\bigskip
\leavevmode

\centerline{\epsfbox{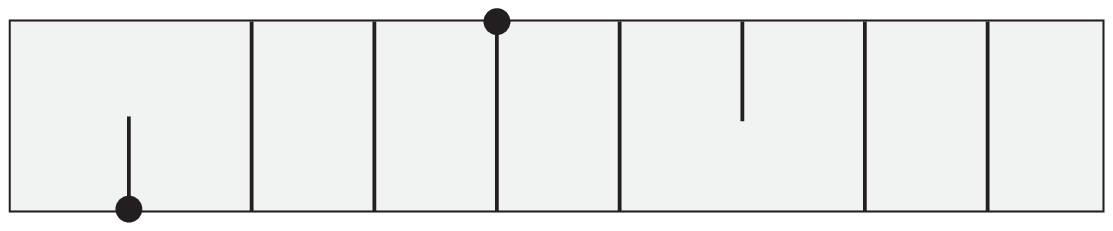}}
\bigskip
\centerline{Figure 3.1}
\bigskip

The parameters of $K(w,b,t)$ appear in Figure 3.1 as follows.  The
winding number $w$ indicates the number of points of $\a$ on the top
or bottom of $R$, labeled successively by $v_1, ..., v_w$ and $v'_1,
..., v'_w$, respectively.  The number $t$ is such that on the torus
$T$ a point $v'_i$ at the bottom is identified to $v_{i+t}$ on the
top, where the subscripts are integers mod $w$.  The bridge number $b$
in $K(w,b,t)$ indicates the number of arcs between the two ending arcs
$\a'$ and $\a''$.  Note that $\b$ is isotopic rel $\bdd$ to an arc in
$R$ with interior intersecting $\a$ at exactly $b$ points in the same
direction.

Now let $\r$ be a $(p,q)$ curve on $T$ such that $M_K(\r)$ is a solid
torus.  By Corollary 2.6, $\b$ is $\r$-trivial, so up to isotopy we
may assume that $\r$ is disjoint from $\a$, and intersects the
meridian disk $D$ above at $p$ points.  Thus $\r \cap A$ is a set of
$p$ vertical arcs.  Cutting $A$ along a component of $\r \cap A$, we
obtain a rectangle $R$ on which $\r$ becomes a set of $p+1$ vertical
arcs $\r_0, ..., \r_p$, with $\r_p$ identified to $\r_0$ on $A$.
These arcs cut $R$ further into a set of rectangles $R_0, ...,
R_{p-1}$.  Since $\a$ is disjoint from $\r$, each arc $\a_i$ lies in
one of the $R_j$.  We may arrange so that the initial point of $\a$
lies in $R_0$, as shown in Figure 3.2, where $(p,q) = (3, 1)$.  Note
that the bottom of $R_i$ is identified to the top of $R_{i+q}$ on $T$
in such a way that the endpoints of the arcs match each other.

\bigskip
\leavevmode

\centerline{\epsfbox{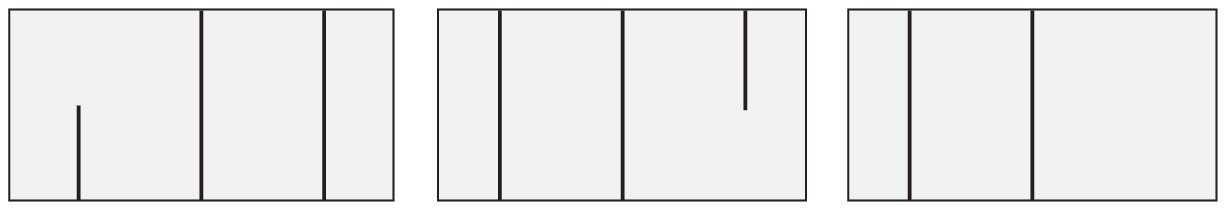}}
\bigskip
\centerline{Figure 3.2}
\bigskip

\definition{Definition 3.1} A 5-tuple $(p,q,k,x,\e)$ of integers is
{\it allowable\/} if 

(1) $p>q > 0$, and $\gcd(p,q) = 1$;

(2) $k \geq 0$;

(3) $p > x \geq 0$;

(4) $\e = 1$  if $k=0$, and $\e \in \{1,-1\}$ if $k>0$.
\enddefinition

Fix a meridian disk $D$ of $V$.  Assume that $\a$ is an arc on $T$
disjoint from a $(p,q)$ curve $\r$ on $T$, such that $\bdd \a \cap
\bdd D = \emptyset$, $\r$ intersects $\bdd D$ minimally up to isotopy,
and $\a$ intersects $\bdd D$ minimally up to isotopy rel $\r \cup \bdd
\a$.  Define an allowable 5-tuple $(p,q,k(\a), x(\a), \e(\a))$ as
follows.

Let $R_i$ be the rectangles defined above.  Up to relabeling we may
assume that $\a(0) \in R_0$.  Let $k(\a)$ be the number of endpoints
of $\a \cap R$ which lie on the top of $R_0$ (hence $\a \cap R_0$ has
$k+1$ components, including $\a'$.)  Let $x = x(\a)$ be the number
such that the arc $\a''$ is contained in $R_x$.

An arc $\a_i$ is a {\it left arc\/} (resp.\ {\it right arc}) if it is
adjacent to the left (resp.\ right) edge of the rectangle $R_j$
containing it.  Define $\e(\a) = 1$ if $\a'$ is a left arc, and
$\e(\a) = -1$ otherwise.  When $k=0$, $\a$ is both a left arc and a
right arc, in which case by definition we have $\e(\a) = 1$.  Note
that one of the $\a'$ and $\a''$ is a left arc, and the other is a
right arc.  For if this were false, assuming that neither of $\a'$ or
$\a''$ is a left arc, say, then the union of the left arcs in $R$
would form a closed circle component of $\a$ on $T$, which would
contradict the fact that $\a$ is an arc on $T$.

We have thus associated an allowable 5-tuple $(p,q, k(\a), x(\a),
\e(\a))$ to each arc $\a$ as above.  Conversely, given a 5-tuple
$(p,q,k,x,\e)$ one can construct such an arc $\a = \a(p,q,k,x,\e)$
with $(p,q,k(\a),x(\a), \e(\a)) = (p,q,k,x,\e)$, which is unique up to
homeomorphism of $(V, \r \cup \bdd D)$.  Let $K= K(p,q,k,x,\e)$ be the
knot isotopic to $\a(p,q,k,x,\e) \cup \b$, where $\b$ is a trivial arc
in $V - D$.  Note that $K(p,q,k,x,\e)$ is a well-defined 1-bridge braid
$K(w,b,t)$ for any allowable 5-tuple.  The parameters $(w,b,t)$ will
be calculated explicitly in Lemma 3.5.  Note that two different
allowable 5-tuples may give rise to the same 1-bridge braid.

The following theorem classifies, for a fixed pair $(p,q)$, all
1-bridge braids $K$ such that $M_K(\r)$ is a solid torus.

\proclaim{Theorem 3.2} Let $K$ be a 1-bridge knot in $V$, and let $p,
q$ be coprime integers such that $0 < q < p$.  Then $M_K(p/q)$ is a
solid torus if and only if $K = K(p,q,k,x,\e)$ for some allowable
5-tuple $(p,q,k,x,\e)$ such that $x = 1$, $q$, $p-q$ or $p-1$.
\endproclaim

\proof Let $\r$ be a $(p,q)$ curve on $T = \bdd V$.  Fix a meridian
disk $D$ of $V$ which intersects $\r$ at $p$ points.  If $M_K(\r)$ is
a solid torus then by Corollary 2.6 the knot $K$ is a 1-bridge braid
$\a\cup \b$, such that $\b$ is $\r$ trivial, and the jumping number
$u(\b, \r)$ is 1 or $q$.  The first condition implies that up to
isotopy one may assume that $\a = \a(p,q,k,x,\e)$ for an allowable
5-tuple $(p,q,k,x,\e)$, and $\b$ is a $\r$-trivial arc disjoint from
$D$.  The second condition implies that $x = 1$, $q$, $p-q$, or $p-1$.
Conversely, if $K = K(p,q,k,x,\e)$ for some allowable 5-tuple
$(p,q,k,x,\e)$ and $x \in \{1, q, p-q, p-1\}$ then one can show that
(1) and (2) in Corollary 2.6 hold, hence $M_K(\r)$ is a solid torus.
\endproof

\proclaim{Corollary 3.3} Given $p>q>0$ with $p\geq 3$, there are
infinitely many 1-bridge braids $K$ in $V$ such that $M_K(p/q)$ is a
solid torus.  \qed \endproclaim

Given a 1-bridge braid $K=K(w,b,t)$, we would like to find all $(p,q)$
such that $M_K(p/q)$ is a solid torus.  By Theorem 3.2 it suffices to
find all allowable 5-tuples $(p,q,k,x,\e)$ such that $K(w,b,t) =
K(p,q,k,x,\e)$, and $x=1,q,p-q$ or $p-1$.  We need to find the
relationship between the parameters $(w,b,t)$ and $(p,q,k,x,\e)$.

Denote by $Z_p$ the set of integers $\{0, 1, ..., p-1\}$, and by
$[a]\in Z_p$ the mod $p$ residue of an integer $a$.  (We will also use
$[x]$ later to denote the integer part of a real number $x$.  It will
be clear from the context what the notation stands for.)  For any
$x \in Z_p$, let $\bar q_x \in Z_p$ be such that $q\bar q_x \equiv x$
mod $p$.  Define $S_x$ as the subset of $Z_p$ given by
$$S_x = \{ [q], [2q], ..., [\bar q_x q] \}.$$ 
Define $\varphi(x,y)$ to be the number of elements $z \in S_x$ which
lie in the open interval $(0,y)$.

Let $\bar q, \bar p$ be the numbers in $Z_p$ such that $\bar q q =
\bar p p + 1$.  Since $\gcd (p,q) = 1$, $\bar q$ and $\bar p$ are well
defined.  The following lemma calculates the value of $\varphi(x,s)$
for certain values of $x$ and $s$.

\proclaim{Lemma 3.4} With the above notation, the values of
$\varphi(x,x)$ and $\varphi(x,q)$ for $x \in \{1, q, p-q, p-1\}$ are
given as follows.

(1)  $\varphi(1,1) = 0$; $\qquad \varphi(1,q) = \bar p$;

(2)  $\varphi(q,q) = 0$;

(3)  $\varphi(p-q, p-q) = p-q -1$; $\qquad \varphi(p-q, q) = q-1$;

(4)  $\varphi(p-1, p-1) = p-\bar q - 1$; $\qquad \varphi(p-1, q) = q - \bar p
     - 1$.
\endproclaim

\proof (1) In this case we have $x = 1$.  By definition $\varphi(1,1)
= 0$.  We need to calculate $\varphi(1,q)$.  The formula is obvious
when $q=1$, so we assume $q>1$.  Note that $\bar q_x = \bar q$.

By definition $\varphi(1,q)$ is the number of points in the set
$\{[q], [2q], ..., [\bar q q] \}$ which are contained in the interval
$(0,q)$.  Since $\gcd(p,q)=1$ and $\bar q < p$, $[iq] \neq 0$ for
$i=1,...,\bar q$.  Hence we need only find the number of $i \in
[1,\bar q]$ such that $[iq] \in [0,q)$.  Note that $[iq] \in [0, q)$
if and only if $iq \in [jp, jp+q)$ for some $j$.

The set of integers $\{q, 2q, ..., \bar q q \}$ distributes evenly on
the interval $[0,\bar q q]$, with distance $q$ between adjacent
numbers.  For any $j=1, 2, ..., \bar p -1$, the interval $[jp, jp+q)$
is contained in $[0,\bar q q)$, hence it contains exactly one $iq$.
Also, the interval $[\bar p p, \bar p p+q)$ contains the point $\bar q
q$ because $\bar q q = \bar p p + 1$ and $q>1$.  When $j>\bar p$, the
interval $[jp, jp+q)$ is disjoint from $[0, \bar q q)$, hence it
contains no $iq$ in the above set.  It follows that exactly $\bar p$
of the numbers in $S_1 = \{ [q], ..., [\bar q q]\}$ are in the
interval $(0,q)$.  Hence $\varphi(1,q) =\bar p$.

(2) When $x=q$, we have $1\cdot q = x$, so $\bar q_x = 1$.  Thus $S_x$
has a single element $[q]$.  Since $[q] = q \notin (0,q)$, we have
$\varphi(x,q) = \varphi(q,q) = 0$.

(3) When $x=p-q$, we have $\bar q_x = p-1$.  Hence $S_x = \{[q], [2q],
..., [\bar q_xq]\}$ contains all integers mod $p$ except $0$.  Thus all
the numbers in $(0,q)$ are in $S_x$, so $\varphi(p-q, q) = q-1$.
Similarly we have $\varphi(p-q, p-q) = p-q-1$.

(4) The calculation for $x=p-1$ is similar to that for $x=1$.  We have
$(p-\bar q) q \equiv p-1$ mod $p$, hence $\bar q_x = p-\bar q$.  From
the equation $(p-\bar q) q = (q-\bar p - 1) p + (p- 1)$, we see that
that $q-\bar p -1$ of the numbers in $S_x = \{[q], ..., [(p-\bar q) q]
\}$ are in the open interval $(0,q)$.  Hence we have $\varphi(p-1,q) =
q-\bar p - 1 $.  Note that all but one of the numbers $\{[q], ...,
[(p-\bar q) q] \}$ are in the interval $(0,p-1)$, (the last number in
the set is $p-1$), so we have $\varphi(p-1, p-1) = p - \bar q - 1$.
\endproof

\proclaim{Lemma 3.5} Let $(p,q,k,x,\e)$ be an allowable 5-tuple.  Then
$K(p,q,k,x,\e) = K(w,b,t)$, where 
$$\gather  w = kp + \bar q_x \\
   t = kq + \varphi(x,q) \\
   b = k(x+\e) + \varphi (x,x)
\endgather
$$
\endproclaim

\proof Recall that $K(p,q,k,x,\e) = \a(p,q,k,x,\e) \cup \b$.
The arc $\a = \a(p,q,k,x,\e)$ is homotopic to $\hat \a \cdot \tilde
\a$, where $\hat \a$ consists of $k$ loops parallel to $\r$, while
$\tilde \a$ is an arc connecting the two endpoints of $\a$ and is
disjoint from the top of the rectangle $R_0$ in Figure 3.2.  The arc
$\tilde \a$ intersects each $R_i$ in either 1 or 0 arc.  By definition
$\bar q_x
\in Z_p$, and $q\bar q_x \equiv x$ mod $p$.  Thus when traveling along
$\tilde \a$ one enters $R_0, R_q, R_{2q}, ..., R_{\bar q_x q}$
successively, where the subscripts are integers mod $p$.  Note that
$\tilde \a \cap R_0$ does not have a vertex on the top of $R_0$,
therefore the wrapping number $w = kp + \bar q_x$.

Recall that $t$ is the number such that the $i$-th point of $\a$ at
the bottom of $R$ is glued to the ($t+i$)-th point of $\a$ on the top
of $R$.  Since the bottom of $R_0$ is glued to the top of $R_q$, the
first point of $\a$ at the bottom of $R_0$ is glued to the first point
of $\a$ on the top of $R_q$, hence $t$ equals the number of points of
$\a$ on the top of $R' = R_0 \cup ... \cup R_{q-1}$.  As above, the
number of points of $\tilde \a$ on the top of $R'$ is $\varphi(x, q)$,
while the number of points of $\hat \a$ on the top of each $R_i$ is
$k$, hence the second equation follows.

When $\e=1$, the bridge width $b$ equal the number of arcs of $\a$ on
$R_0\cup ...\cup R_x$, not counting the two arcs containing an
endpoint of $\a$.  Therefore we have $b = k(x+1) + \varphi(x,x)$.
When $\e=-1$, the bridge does not pass the edges in $R_0$ and $R_x$,
hence $b = k(x-1) + \varphi(x,x)$.
\endproof

The following theorem gives a simple method to calculate, for a given
1-bridge braid $K$, all $(p,q)$ such that $M_K(\r)$ is a solid torus,
where $\r$ is a $(p,q)$ curve.  There are four possible solutions for
each $K$; in each case $(p,q)$ can be calculated from $w$ and $t$,
and it has to satisfy an extra condition in terms of $b$.  

\proclaim{Theorem 3.6} Let $K = K(w,b,t)$ be a 1-bridge braid in a
solid torus $V$ such that $1\leq b \leq w-2$, and $0 < t < w$.  Denote
by $[y]$ the integer part of a real number $y$.  Then $M_K(p/q)$ is a
solid torus if and only if $p,q$ satisfy one of the following
conditions.

(1)  $qw - pt = 1$, $p,q \in Z_w$, and  $b = 2[w/p]$.

(2) $k = \gcd (w-1, t)$, $p = (w-1)/k$, $q = t/k$, and $b = k(q+\e)$ for
some $\e=\pm 1$.

(3) $k = \gcd (w+1, t+1) - 1$, $p = (w+1)/(k+1)$, $q = (t+1)/(k+1)$, and $b
= k(p-q+\e) + (p-q -1)$ for some $\e=\pm 1$.

(4) $p(t+1) - qw = 1$, $p,q \in Z_w$, and $b= [w/p](p-2) +
(p-\bar q-1)$.  
\endproclaim

\proof Let $M = M_K$ be the exterior of $K$.  By Theorem 3.2, $M(\r)$
is a solid torus if and only if $K(w,b,t) = K(p,q,k,x,\e)$ for some
allowable 5-tuple $(p,q,k,x,\e)$, and $x = 1,q, p-q$ or $p-1$.  We
need to show that these correspond to (1) -- (4) in the theorem.

\medskip

CASE 1:  $x=1$.  In this case by Lemma 3.4 the three equations in
Lemma 3.5 become
$$\align  w &= kp + \bar q_1 = kp + \bar q\\
   t &= kq + \varphi(1,q) = kq + \bar p\\
   b &= k(1+\e) + \varphi (1,1) = k(1+\e)
\endalign
$$
If $K(w,b,t) = K(p,q,k,x,\e)$ and $x = 1$, then since $b>0$, the
third equation above gives $\e=1$, so $b=2k=2[w/p]$.  We have
$qw - pt = q(kp + \bar q) - p(kq + \bar p) = q \bar q - p \bar p =
1$.  Therefore (1) holds.

Conversely, if $p,q \in Z_w$ satisfy the equation $qw - pt =
1$, let $w = k' p + p'$, and $t = k'' q + q'$, where $0\leq p'< p$
and $0 \leq q' < q$.  Then the equation $qw - pt = 1$ becomes
$$ 1 = qw - pt = q(k' p + p') - p(k'' q + q') = (k' - k'') pq + (qp' -
pq').$$ 
Since $p,q \in Z_w$, this implies that $k'=k''$, and $qp' - pq' = 1$;
hence $p' = \bar q$ and $q' = \bar p$.  Hence if we define $k = k' =
k''$ and $x=\e=1$ then $(p,q,k,x,\e)$ is an allowable 5-tuple with
$x=1$.  By definition we have $w = kp + \bar q$ and $t = kq + \bar p$.
The equation $b = k(1 + \e)$ follows from the extra condition $b = 2
[w/p] = 2k$ in (1).  Therefore $K(w,b,t) = K(p,q,k,x,\e)$.

The arguments for the other cases are similar.  We only show the
calculations below.

\medskip

CASE 2:  $x=q$.  By Lemma 3.4 the equations in Lemma 3.5 become
$$\align  w &= kp + \bar q_x = kp + 1\\
   t &= kq + \varphi(x,q) = kq \\
   b &= k(x+\e) + \varphi (x,x) = k(q+\e)
\endalign
$$
We can solve the first two equations for $k,p,q$ to get
$k=\gcd(w-1,t)$, $p=(w-1)/k$, and $q=t/k$.  The third equation gives
the extra condition that must be satisfied, i.e., $b=k(q+\e)$ for some
$\e=\pm 1$.

\medskip

CASE 3:  $x=p-q$.  In this case $\bar q_x = p-1$.  We have
$$\align  w &= kp + \bar q_x = kp + (p-1)\\
   t &= kq + \varphi(x,q) = kq + (q-1)\\
   b &= k(x+\e) + \varphi (x,x) = k(p-q+\e)) + (p-q-1)
\endalign
$$
The first two equations give $k=\gcd(w-1,t)-1$, $p=(w-1)/(k+1)$, and
$q=t/(k+1)$.  The third equation gives the extra condition that 
$b=k(p-q+\e) + (p-q+1)$ for some $\e=\pm 1$.

\medskip

CASE 4:  $x = p-1$.  
In this case $\bar q_x = p-\bar q$.  We have
$$\align  w &= kp + \bar q_x = kp + (p-\bar q)\\
   t &= kq + \varphi(x,q) = kq + (q-\bar p - 1)\\
   b &= k(x+\e) + \varphi (x,x) = k(p-1+\e)) + (p-\bar q-1)
\endalign
$$
As in Case 1, one can show that the first two equations are
equivalent to the condition that $p(t+1) - qw = 1$ for some $p,q\in
Z_w$.  As usual, we have $k=[w/p]$, hence the last equation is
equivalent to $b= [w/p](p-1+\e)) + (p-\bar q-1)$.  Since $b \neq w-1$,
we must have $\e = -1$.  \endproof

\proclaim{Corollary 3.7} Let $K = K(w,b,t)$ be a 1-bridge braid in a
solid torus $V$ such that $1\leq b \leq w-2$, and $0 < t < w$.  If
$M_K (p/q)$ is a solid torus, then $w+1 \geq p > q > 0$.  \qed
\endproclaim

\head 4. Computational results \endhead

A computer program can be written using Theorems 3.2 and 3.6 to
calculate a list, for $w$ up to a given value, of all $K(w,b,t)$ and
$(p,q)$ such that $(p,q)$ filling on the outer torus of $M_K$ produces
a solid torus.  It can also calculate all such $(p,q)$ for any given
$K(w,b,t)$.  The following are some results from such a calculation.

(1) The exterior of the knot $K(7,2,4)$ admits three such fillings,
    with slopes $(3,2)$, $(5,3)$ and $(8,5)$.  Since the dual knot
    after such a filling is a hyperbolic knot in a solid torus, by a
    theorem of Berge [Be], up to (possibly orientation reversing)
    homeomorphism of $(V, K)$ this is the only 1-bridge braid such
    that $M_K$ admits three solid torus fillings on $T$.  Note that
    $K(7,4,2)$ is equivalent to $K(7,2,4)$ by an orientation reversing
    map of $V$.

(2) The knot $K(8,3,6)$ is the first one (in lexicographic order of
    parameters) that admits no solid torus filling on $T$.

(3) As $w$ increases, the percentage of knots $K(w,b,t)$ admitting
    solid torus fillings on $T$ becomes smaller, as expected, and most
    of them only admit one such filling, so its dual knot does not
    admit nontrivial surgery which produces a solid torus.  There are
    72 knots with $w\leq 10$, 60 of which admit a total of 86 solid
    torus fillings on $T$.  There are exactly 6000 knots $K(w,b,t)$
    with $w\leq 40$, 2380 of which admit a total of 2692 such
    fillings.

One can easily show that there is an orientation reversing map of $V$
sending $K=K(w,b,t)$ to $\bar K = K(w,w-b-1, t-b-1)$, which sends a
$(p,q)$ curve on $T$ to a $(p,p-q)$ curve, so $M_K(p/q)$ is a solid
torus if and only if $M_{\bar K}(p/(p-q))$ is.  Thus up to orientation
reversing homeomorphism of $V$ we may assume that $b < w/2$, and if
$b=(w-1)/2$ then $t<w/2$.  There are 36 such knots for $w \leq 10$.
The following list shows all such knots and the Dehn filling slopes on
$T$ such that $M_K(p/q)$ is a solid torus.  As an example, the tuple
$(6,2,3; 5/3, 7/4)$ indicates that the knot $K(6,2,3)$ admits two such
Dehn fillings, with slopes $5/3$ and $7/4$ respectively.

\bigskip

\centerline{TABLE 1.  $(w,b,t; p_1/q_1, p_2/q_2, p_3/q_3)$}

$$ \matrix
(4, 1, 2; 3/2, 5/3)  &
(5, 2, 1; 3/1, 4/1)   &
(6, 1, 2; 5/2)   \\
(6, 1, 4; 7/5)   &
(6, 2, 1; 5/1)   &
(6, 2, 3; 5/3, 7/4) \\
(7, 2, 1; 6/1)   &
(7, 2, 4; 3/2, 5/3, 8/5)   &
(7, 2, 5; 4/3)   \\
(8, 1, 2; 7/2)   &
(8, 1, 6; 9/7)   &
(8, 2, 1; 7/1)   \\
(8, 2, 3; 5/2, 7/3)   &
(8, 3, 2; 3/1, 7/2)   &
(8, 3, 4; 7/4, 9/5)   \\
(8, 3, 6; - )   &
(9, 2, 1; 8/1)   &
(9, 2, 3; 8/3)   \\
(9, 2, 5; 5/3, 7/4)   &
(9, 2, 6; 10/7)   &
(9, 2, 7; 5/4)   \\
(9, 4, 1; -)   &
(9, 4, 2; 4/1)   &
(9, 4, 3; 5/2, 8/3)   \\
(10, 1, 2; 9/2)   &
(10, 1, 4; -)   &
(10, 1, 6; -)   \\
(10, 1, 8; 11/9)   &
(10, 2, 1; 9/1)   &
(10, 2, 7; 7/5, 11/8)   \\
(10, 3, 2; 9/2)   &
(10, 3, 4; 9/4)   &
(10, 3, 6; 3/2, 11/7)   \\
(10, 3, 8; -)   &
(10, 4, 1; -)   &
(10, 4, 5; 9/5, 11/6)   
\endmatrix
$$

{\it Acknowledgement.\/} I would like to thank the referee for his/her
careful reading, and for some very helpful comments.

\enddocument